# LOCAL-GLOBAL PROBLEM FOR DRINFELD MODULES

G.J. VAN DER HEIDEN

ABSTRACT. In this paper we study the group
$$S(a,K) := \ker\left(E(K)/aE(K) \longrightarrow \prod_\nu E(K_\nu)/aE(K_\nu)\right)$$
in two cases. In the first case $E$ is a Drinfeld module, $K$ a function field and $(a)$ a principal prime ideal in the ring $A$ of regular functions which occurs in the definition of a Drinfeld module. In the second case $E$ is an elliptic curve over a number field $K$ and $a$ is a prime number. We will prove that $S(a,K)$ is trivial in the elliptic curve case. In the Drinfeld module case we will compute $S(a,K)$ when the Drinfeld module has rank 2. We will show in what cases this group is trivial and that it can become arbitrarily large.

## 1. THE DRINFELD MODULE CASE

Let $X$ be a projective, smooth, absolutely irreducible curve over $\mathbb{F}_q$, with $\operatorname{char}(\mathbb{F}_q) = p$, and let $\infty \in X$ be some fixed closed point on $X$. $\mathbb{F}_q(X)$ will denote the function field of $X$ and $A$ will denote the ring of functions in $\mathbb{F}_q(X)$ which are regular outside $\infty$.

Let $K$ be some separable, finite extension of $\mathbb{F}_q(X)$, which is also an $A$-algebra via the natural embedding $\gamma : A \longrightarrow K$ and let $K^s$ be the separable closure of $K$ inside some algebraic closure $K^{\operatorname{alg}}$ of $K$. Let $K\{\tau\}$ be the skew polynomial ring consisting of elements $\sum_i k_i \tau^i$, $k_i \in K$. Multiplication in $K\{\tau\}$ is given by the rule $\tau k = k^q \tau$ for all $k \in K$.

Let $\varphi$ be a *Drinfeld module over $K$ of rank $r$*, i.e. $\varphi$ is a ring homomorphism
$$\varphi : A \longrightarrow \operatorname{End}_{\mathbb{F}_q}(\mathbb{G}_{a,K}) \cong K\{\tau\},$$
such that for all $a \in A$
$$\varphi_a = \sum_{i=0}^{r \deg(a)} k_i \tau^i, \quad \text{with } k_{r\deg(a)} \in K^*, k_0 = \gamma(a).$$

Note that we write $\varphi_a$ instead of $\varphi(a)$.

For any field extension $L$ of $K$, $\varphi$ induces an $A$-module structure on $L$ via the embedding $K \subset L$, by
$$\varphi_a(l) = \sum_i k_i \tau^i(l) = \sum_i k_i l^{q^i} \quad \text{for all } l \in L.$$

*Date*: May 30, 2002.

$E(L)$ denotes $L$ with this $A$-module structure.

We write $E[a](L)$ for $\ker(\varphi_a)(L)$ and moreover $aQ$ for $\varphi_a(Q)$, with $Q \in E(L)$.

Let $(a) \subset A$ be a principal prime ideal of $A$. The object of study in this paper is the kernel

$$S(a, K) := \ker\left(E(K)/aE(K) \longrightarrow \prod_\nu E(K_\nu)/aE(K_\nu)\right),$$

where the product is taken over all places $\nu$ of $K$.

## 2. The group $S(a, K)$

For any $P \in E(K)$, let

$$K_P := \text{the splitting field of } \varphi_a(Z) - P \in K[Z] \text{ over } K.$$

**Lemma 2.1.** *The finite field extensions $K \subset K_0 \subset K_P$ are Galois.*

*Proof.* Because $\frac{d}{dZ}(\varphi_a(Z) - P) = \gamma(a) \neq 0$, the polynomial $\varphi_a(Z) - P$ is separable and hence $K_P$ is a finite Galois extension of $K$. In particular $K_0$ is a Galois extension of $K$. Moreover for any $P$, let $Q \in K_P$ be a zero of $\varphi_a(Z) - P$, then $K_P = K_0(Q)$, because $Q + \ker(\varphi_a(Z))$ is the set of all zeroes of $\varphi_a(Z) - P$. Hence the extension $K_P \supset K_0$ is also Galois. □

Note that $\mathbb{F}_q(X)$ is a function field, i.e., a finite separable extension of $\mathbb{F}_q(t)$, for some transcendental element $t \in \mathbb{F}_q(X)$. For function fields, as well as for number fields, we have *Chebotarev's density theorem*, cf. [Jar82]. The following lemma is a consequence of this theorem.

**Lemma 2.2.** *Let $K$ be a function field (resp. a number field) and let $L$ be a finite separable extension of $K$. If for all places $\nu$ of $K$ there exists a place $\omega$ of $L$ of degree 1, then $K = L$.*

*Proof.* Let $M$ be the normal closure of $L/K$, then both $M/L$ and $M/K$ are finite Galois extensions. Let $H = \text{Gal}(M/L)$ and $G = \text{Gal}(M/K)$. By Chebotarev, every $\sigma \in G$ is the Frobenius of some place $\mu$ of $M$ lying above some place $\nu$ of $K$. This implies that $\sigma \mod \nu$ generates the Galois group $\text{Gal}(k_\mu/k_\nu)$, where $k_\mu$ and $k_\nu$ denote the residue fields at $\mu$ and $\nu$ respectively. Because both $M/K$ and $M/L$ are Galois, there is a $\tau \in G$ such that the conjugate $\mu' = \tau(\mu)$ of $\mu$ lies above a place $\omega$ of $L$, which has degree 1 over $\nu$. In particular we see that $\tau\sigma\tau^{-1}$ generates $\text{Gal}(k_{\mu'}/k_\nu) = \text{Gal}(k_{\mu'}/k_\omega)$, where the equality follows from the fact that $\deg(\omega/\nu) = 1$. And thus we see that $\tau\sigma\tau^{-1} \in H$, so $\sigma \in \tau^{-1}H\tau$.

We conclude that

$$G = \bigcup_{\tau \in G} \tau H \tau^{-1} = \bigcup_{\tau H \in G/H} \tau H \tau^{-1}.$$

Note that $1 \in \tau H \tau^{-1}$ for all $\tau \in G$. On the other hand, $G$ equals the union of all distinct cosets $\tau H$, which is a disjoint union. By comparing the number of elements one sees that this is only possible if $H = G$, hence $K = M^G = M^H = L$. □

**Proposition 2.3.** *For every class $[P] = P + aE(K) \in S(a, K)$ we have $K_P = K_0$. In particular $S(a, K_0) = 0$.*

*Proof.* First note that for every $Q \in aE(K)$, we have $K_Q = K_0$, hence the extension $K_P$ only depends on the class $[P] = P + aE(K)$.

Let now $P \in S(a, K)$ and let $\nu$ be a place of $K$. Then $K \subset K_0 \subset K_P$ and correspondingly we have places $\nu, \nu_0$ and $\nu_P$, with $\nu_0$ a place of $K_0$ lying above $\nu$ and $\nu_P$ a place of $K_P$ lying above $\nu_0$.

Because $P \in S(a, K)$, there exists a solution of $\varphi_a(Z) - P = 0$ in $K_\nu$, hence all solutions of this equation lie in $(K_0)_{\nu_0}$. This means that $K_0 \subset K_P \subset (K_0)_{\nu_0}$. It follows that $(K_0)_{\nu_0} = (K_P)_{\nu_P}$ and in particular $\deg(\nu_P/\nu_0) = 1$. From lemma's (2.1) and (2.2), it follows that $K_P = K_0$.

If $[P] \in S(a, K_0)$, then $(K_0)_P = (K_0)_0 = K_0$, hence $P \in aE(K_0)$, thus $[P] = 0$. □

We write thoughout this paper $G_{K_0} = \text{Gal}(K_0/K)$.

**Proposition 2.4.** *We have*

$$S(a, K) \cong \bigcap_\omega \ker(\text{Res}_\omega).$$

*Here the intersection is taken over all places $\omega$ of $K_0$ and $\text{Res}_\omega$ is the restriction map*

$$\text{Res}_\omega : H^1(G_{K_0}, E[a](K_0)) \longrightarrow H^1(D_\omega, E[a]((K_0)_\omega)),$$

*where $D_\omega$ denotes the decomposition group at $\omega$.*

*Proof.* Consider the following diagram:

$$
\begin{array}{ccccccc}
 & & 0 & & 0 & & \\
 & & \downarrow & & \downarrow & & \\
0 & \longrightarrow & C & \longrightarrow & S(a, K) & \longrightarrow & S(a, K_0) = 0 \\
 & & \downarrow & & \downarrow & & \downarrow \\
0 & \longrightarrow & \Phi & \longrightarrow & E(K)/aE(K) & \longrightarrow & E(K_0)/aE(K_0) \\
 & & \downarrow & & \downarrow & & \downarrow \\
0 & \longrightarrow & B & \longrightarrow & \prod_\nu E(K_\nu)/aE(K_\nu) & \longrightarrow & \prod_\omega E(K_{0,\omega})/aE(K_{0,\omega}).
\end{array}
$$

Clearly, the second and third row are exact as well as all columns. From this it follows that also the first row is exact. Hence we see that $S(a, K) \cong C$.

To determine the groups $C$ and $B$, we use some Galois cohomology. Starting from the exact sequence

$$0 \longrightarrow E[a](K^s) \longrightarrow E(K^s) \xrightarrow{\varphi_a} E(K^s) \longrightarrow 0,$$

we deduce that

$$E(K)/aE(K) \hookrightarrow H^1(\text{Gal}(K^s/K), E[a](K^s)).$$

By additive Hilbert 90 the cokernel of this map is $H^1(\mathrm{Gal}(K^s/K), E(K^s)) = 0$, because $E(K^s)$ is here just $(K^s)^+$, hence

$$E(K)/aE(K) \cong H^1(\mathrm{Gal}(K^s/K), E[a](K^s)).$$

Similarly we deduce that

$$E(K_0)/aE(K_0) \cong H^1(\mathrm{Gal}(K^s/K_0), E[a](K^s)).$$

In particular this means that

$$\Phi \cong H^1(G_{K_0}, E[a](K_0)).$$

Applying the same arguments to $K_\nu$ and $K_{0,\omega}$, for some place $\nu$ of $K$ and $\omega$ of $K_0$ lying above $\nu$, we deduce

$$\ker\left(E(K_\nu)/aE(K_\nu) \longrightarrow E(K_{0,\omega})/aE(K_{0,\omega})\right) \cong H^1(D_\omega, E[a](K_{0,\omega})),$$

where $D_\omega$ denotes the decomposition group at $\omega$. This isomorphism implies that

$$B \cong \prod_\nu \bigcap_{\omega|\nu} H^1(D_\omega, E[a](K_{0,\omega})),$$

where the product runs over all places $\nu$ of $K$. If we furthermore note that the map $\mathrm{Res}_\omega$ only depends on the place $\nu$ underlying $\omega$, it follows that $C$ is the kernel of the map $\prod_\nu \mathrm{Res}_\omega$, with $\mathrm{Res}_\omega$ as in the proposition. □

## 3. The group $H^1(G_{K_0}, E[a](K_0))$

In the following we will write $\mathbb{F} = A/(a)$ and $V = E[a](K_0)$. For all field extensions $L \supset K_0$, we have $V = E[a](L)$. Note that $\mathbb{F}$ is a field extension of $\mathbb{F}_q$, because $(a)$ is prime. It is well-known that $V \cong \mathbb{F}^r$, where $r$ is the rank of $\varphi$. The action of $\sigma \in G_{K_0}$ on elements in $K_0$ commutes with the action of $\varphi_f$ for all $f \in A$ and hence we have a representation

$$G_{K_0} \hookrightarrow \mathrm{Gl}_r(\mathbb{F}),$$

which is an embedding because $K_0$ is given by adjoining the elements of $V$, which are the zeroes of $\varphi_a(Z)$, to $K$.

**Proposition 3.1.** *For every Drinfeld module of rank 1, we have $S(a, K) = 0$.*

*Proof.* Note that $G_{K_0} \hookrightarrow \mathbb{F}^*$ and thus $p \nmid \#G_{K_0}$, but $p$ does divide $\#V = \#\mathbb{F}$, hence $H^1(G_{K_0}, V) = 0$. □

**Proposition 3.2.** *Let $\mathbb{F}$ be a finite field of characteristic $p$ and let $W$ be an $\mathbb{F}$-vectorspace of dimension $r$. If $\mathbb{F} \neq \mathbb{F}_2$, then*

$$H^1(\mathrm{Gl}_r(\mathbb{F}), W) = 0.$$

*If $p \neq 2$, then*

$$H^1(\mathrm{Sl}_r(\mathbb{F}), W) = 0.$$

*Proof.* For the first part, note that if $\mathbb{F} \neq \mathbb{F}_2$, then we may choose $\alpha \in \mathbb{F}^*$, such that $\alpha \neq 1$. Hence $H = \langle \alpha I \rangle$ is a non-trivial normal subgroup of $\mathrm{Gl}_r(\mathbb{F})$. Note that $W^H = 0$. Moreover $H^1(H, W) = 0$, because this group is annihilated by both $\#H$, which is prime to $p$, and $p$. By the exact sequence

$$0 = H^1(\mathrm{Gl}_r(\mathbb{F})/H, W^H) \longrightarrow H^1(\mathrm{Gl}_r(\mathbb{F}), W) \longrightarrow H^1(H, W) = 0,$$

the first statement follows.

For the second part, we let $H$ be the normal subgroup of $\mathrm{Sl}_r(\mathbb{F})$ consisting of the diagonal matrices, whose diagonal entries are $1$ or $-1$. If $r > 1$, then $H$ is not trivial and $W^H = 0$. Moreover $\#H \mid 2^r$, hence $H^1(H, W) = 0$, because it is annihilated by both $2$ and $p \neq 2$. Applying the same sequence as in the proof of the first part gives the result. □

*Remark* 3.3. For rank $r = 2$ the Galois group $G_{K_0}$ is generically $\mathrm{Gl}_2(\mathbb{F})$, cf. [Gar01]. It is conjectured that for arbitrary rank this is also true, i.e. the Galois group is generically $\mathrm{Gl}_r(\mathbb{F})$. Proposition (3.2) states that given this conjecture, $S(a, K)$ is generically $0$.

## 4. The rank 2 case.

From now on we will assume that the rank of the Drinfeld module $\varphi$ is 2. Throughout the rest of this paper $H$ will denote

$$H := G_{K_0} \cap \mathrm{Sl}_2(\mathbb{F}).$$

Note that $H$ is a normal subgroup of $G_{K_0}$ and that $p \nmid [G_{K_0} : H]$, hence

$$H^1(G_{K_0}/H, V^H) = 0$$

and we see by group cohomology that

$$H^1(G_{K_0}, V) \hookrightarrow H^1(H, V).$$

The classification of subgroups $\mathrm{Sl}_2(\mathbb{F}_q)$ from [Suz82] shows that $H$ is one of the following.

(1) $p \nmid \#H$.
(2) $D_{2n}$; in this case $p = 2$ and $n$ is odd.
(3) $A_5$; in this case $p = 3$.
(4) $\mathrm{Sl}_2(\mathbb{F}_{p^k})$, where $\mathbb{F}_{p^k} \subset \mathbb{F}_q$.
(5) $\langle \mathrm{Sl}_2(\mathbb{F}_{p^k}), \begin{pmatrix} \lambda & 0 \\ 0 & \lambda^{-1} \end{pmatrix} \rangle$, where $\lambda^2$ generates $\mathbb{F}_{p^k}$, but $\lambda \notin \mathbb{F}_{p^k}$.
(6) $H$ is a Borel group, i.e., $H$ has a normal abelian $p$-Sylow subgroup $Q$, such that $H/Q$ is cyclic of order dividing $\#\mathbb{F}^*$.

In the following proposition we deal with most of the subgroups in the classification.

**Theorem 4.1.** *If $H$ is of type (1) or (2) or if $p > 2$ and $2 \mid \#H$, then*

$$H^1(H, V) = 0.$$

*Consequently, in all these cases $S(a, K) = 0$.*

*Proof.* We consider the different types of $H$:

*Type (1).* $H^1(H, V)$ is annihilated by both $p$ and $\#H$, hence
$$H^1(H, V) = 0.$$

*Type (2).* In this case $p = 2$. We consider the following exact sequence, where $x$ is one of the generators of order 2 of $D_{2n}$
$$H^1(D_{2n}, V) \xrightarrow{\text{Res}} H^1(\langle x \rangle, V) \xrightarrow{\text{Cor}} H^1(D_{2n}, V).$$
Now by [Ser79] Cor $\circ$ Res $= n$. Because $x$ is of order 2 and $p = 2$, we may assume that
$$x = \begin{pmatrix} 1 & c \\ 0 & 1 \end{pmatrix}.$$
Proposition (4.3) implies that
$$H^1(\langle x \rangle, V) = 0,$$
hence $H^1(H, V)$ is annihilated by both 2 and an odd number $n$.

*Type $p > 2$ and $2 \mid \#H$.* This implies that $H$ contains the non-trivial normal subgroup $\langle \pm 1 \rangle$. Now the exact sequence
$$0 = H^1(H/\langle \pm 1 \rangle, V^{\langle \pm 1 \rangle}) \longrightarrow H^1(H, V) \longrightarrow H^1(\langle \pm 1 \rangle, V) = 0,$$
gives that $H^1(H, V) = 0$. $\square$

*Remark 4.2.* The only cases of the classification which are not covered by this theorem are the cases when $H$ is of type (4) or (5), when $p = 2$, and $H$ is of type (6) (such that $2 \nmid \#H$).

If $p = 2$ and $H$ is of type (4), then by [CPS75], table 4.5, we obtain that $\dim_{\mathbb{F}} H^1(G_{K_0}, V) = 1$. So if $p = 2$ and $H$ is of type (4) or (5), this might give rise to examples for which $S(a, K)$ is non-trivial. If this is possible, then the construction should be analogous to the construction in case $H$ is of type (6), which we will give in section (5).

*H of type (6).* In the rest of this section we will assume that $H$ is of type (6) and we compute $H^1(H, V)$. Let $Q$ by the $p$-Sylow subgroup of $H$. Clearly
$$H^1(H/Q, V^Q) = 0,$$
because this group is annihilated by both $p$ and $\#(H/Q)$, which is prime to $p$. It follows that
$$H^1(G_{K_0}, V) \hookrightarrow H^1(H, V) \hookrightarrow H^1(Q, V).$$
Let $k \in \mathbb{N}$ such that $p^k = \#Q$, then $Q = \langle \sigma_1, \ldots, \sigma_k \rangle$ and $H = \langle Q, \rho \rangle$, where
$$\sigma_i = \begin{pmatrix} 1 & \lambda_i \\ 0 & 1 \end{pmatrix}, \quad \rho = \begin{pmatrix} \alpha & 0 \\ 0 & \alpha^{-1} \end{pmatrix},$$

such that the $\lambda_i$ are linearly independent over $\mathbb{F}_p$ and $\alpha \in \mathbb{F}^*$ generates $H/Q$. Let $\tau \in Q$ and write $\mathrm{Res}_{\langle \tau \rangle}$ for the residue map

$$\mathrm{Res}_{\langle \tau \rangle} : H^1(H, V) \longrightarrow H^1(\langle \tau \rangle, V).$$

**Proposition 4.3.** *The $\mathbb{F}$-vectorspace $H^1(Q, V)$ has dimension*

$$\dim_{\mathbb{F}} H^1(Q, V) = \begin{cases} \dim_{\mathbb{F}_p} Q & \text{if } p > 2 \\ -1 + \dim_{\mathbb{F}_p} Q & \text{if } p = 2. \end{cases}$$

*If $H = Q$ and $\sigma \in Q$ is not the identity, then*

$$\dim_{\mathbb{F}} \ker(\mathrm{Res}_{\langle \sigma \rangle}) = -1 + \dim_{\mathbb{F}_p} Q.$$

*Proof.* We write $V = \mathbb{F}e_1 + \mathbb{F}e_2$. Note that $V$ is an $\mathbb{F}[Q]$-module. The group ring $\mathbb{F}[Q]$ is isomorphic to the commutative ring $\mathbb{F}[\sigma_1, \ldots, \sigma_k]$ with only relations $\sigma_i^p = 1$. Write $x_i = \sigma_i - 1$, then $\mathbb{F}[Q]$ is the commutative ring $R = \mathbb{F}[x_1, \ldots, x_k]$ subject to the relations $x_i^p = 0$ and these are the only relations.

Note that $\mathbb{F}$ is isomorphic to $R/(x_1, \ldots, x_k)$. To compute $H^1(Q, V)$, we consider the following free resolution of the $R$-module $\mathbb{F}$:

$$R^{k + \frac{1}{2}k(k-1)} \xrightarrow{d_1} R^k \xrightarrow{d_0} R \xrightarrow{d_{-1}} \mathbb{F} \longrightarrow 0,$$

where $R$ linear maps are given as follows:

$$d_{-1} : 1 \mapsto 1 \mod (x_1, \ldots, x_k),$$

write $b_1, \ldots, b_k$ for generators of $R^k$ over $R$, then

$$d_0 : b_i \mapsto x_i,$$

write $c_1, \ldots, c_k, c_{i,j}$ with $1 \leq i < j \leq k$ for the generators of $R^{k + \frac{1}{2}k(k-1)}$, then

$$d_1 : c_i \mapsto x_i^{p-1} b_i, \quad d_1 : c_{i,j} \mapsto x_i b_j - x_j b_i.$$

To see that the given sequence is exact, note that $\ker(d_0)$ is generated by the elements $d_1(c_i), d_1(c_{i,j})$ for all $i, j$, because these exactly describe all relations in de the ring $R$.

From this sequence we arrive at the cocomplex

$$V^{k + \frac{1}{2}k(k-1)} \xleftarrow{d_1} V^k \xleftarrow{d_0} V \longleftarrow 0,$$

with

$$d_0(v) = (x_1 v, \ldots, x_k v),$$

and

$$d_1(v_1, \ldots, v_k) = (x_1^{p-1} v_1, \ldots, x_k^{p-1} v_k, (x_i v_j - x_j v_i)_{i<j}).$$

To compute $\ker(d_1)$ and $\mathrm{im}(d_0)$, note that the action of $R$ on $V$ is given by $x_i e_1 = 0$ and $x_i e_2 = \lambda_i e_1$ for all $i$. From this it follows immediately that $\mathrm{im}(d_0)$ is generated over $\mathbb{F}$ by the vector $(\lambda_1 e_1, \ldots, \lambda_k e_1)$. Hence $\dim_{\mathbb{F}} \mathrm{im}(d_0) = 1$.

To compute $\ker(d_1)$, note that $x_i^{p-1} v = 0$ for all $v$ if $p > 2$. So if $p > 2$, then an element of $V^k$ lies in $\ker(d_1)$ iff $x_i v_j = x_j v_i$. Write $v_i = a_i e_1 + b_i e_2$, with $a_i, b_i \in \mathbb{F}$, then

$$x_i(a_j e_1 + b_j e_2) = \lambda_i b_j e_1 = x_j(a_i e_1 + b_i e_2) = \lambda_j b_i e_1.$$

From this it follows that $\ker(d_1)$ is generated by

$$(e_1, 0, \ldots, 0), \ldots, (0, \ldots, 0, e_1), (\lambda_1 e_2, \ldots, \lambda_k e_2),$$

hence $\dim_{\mathbb{F}} \ker(d_1) = k+1$. So we see that for $p > 2$, that $\dim_{\mathbb{F}} H^1(Q, V) = k$.

If $p = 2$, then elements in $\ker(d_1)$ must satisfy $x_i^{p-1} v_i = x_i v_i = 0$, hence $v_i = a_i e_1$, with $a_i \in \mathbb{F}$. Hence $\ker(d_1)$ is contained in the span of

$$(e_1, 0, \ldots, 0), \ldots, (0, \ldots, 0, e_1).$$

For vectors in this span clearly also the other equations $x_i v_j = x_j v_i$ hold, hence this span equals $\ker(d_1)$ and thus $\dim_{\mathbb{F}} \ker(d_1) = k$. So for $p = 2$, we have $\dim_{\mathbb{F}} H^1(Q, V) = k - 1$.

Clearly, by this computation $\dim_{\mathbb{F}} H^1(\langle \sigma \rangle, V) = 1$ if $p > 2$ and $0$ if $p = 2$. From this the dimension formula for $\ker(\mathrm{Res}_{\langle \sigma \rangle})$ follows. $\square$

**Proposition 4.4.** *Suppose that $H \neq Q$, say $H/Q \cong \langle \alpha \rangle$, with $2 \nmid \mathrm{ord}(\alpha)$. Let $\delta = 1$ if $\mathrm{ord}(\alpha) = 3$ and $p > 2$ and $\delta = 0$ otherwise. Let $l = \dim_{\mathbb{F}_p[\alpha]} \mathbb{F}$. Then*

$$\dim_{\mathbb{F}} H^1(H, V) = \begin{cases} 0 & \text{if } \alpha^{p^j} \neq \alpha^2 \text{ for all } j; \\ l + \delta & \text{otherwise.} \end{cases}$$

*Proof.* First note that we may extend the restriction-inflation sequence as follows (cf. [Ser79]):

$$0 \longrightarrow H^1(H/Q, V^Q) \longrightarrow H^1(H, V) \longrightarrow H^1(Q, V)^{H/Q} \longrightarrow H^2(H/Q, V^Q),$$

which induces an isomorphism

$$H^1(H, V) \cong H^1(Q, V)^{H/Q},$$

because $H^1(H/Q, V^Q) = H^2(H/Q, V^Q) = 0$.

We will now compute the $H/Q$-invariant cocycle classes in $H^1(Q, V)$. We will use the following notation: for a cocycle $\xi : Q \longrightarrow V$ we write

$$\xi(\sigma_\lambda) = \begin{pmatrix} x_\lambda \\ y_\lambda \end{pmatrix}, \quad \text{where} \quad \sigma_\lambda = \begin{pmatrix} 1 & \lambda \\ 0 & 1 \end{pmatrix}.$$

We write $x : \mathbb{F} \longrightarrow \mathbb{F}$ for the map $x : \lambda \mapsto x_\lambda$ and $y : \mathbb{F} \longrightarrow \mathbb{F}$ for the second coordinate, then the cocycle relations and the relations between the elements in $Q$ imply that $x$ and $y$ are determined by the following relations:

(1a) $\quad x(\mu + \lambda) = x(\mu) + x(\lambda) + \lambda \mu y(1)$

(1b) $\quad y(\lambda) = \lambda y(1) \quad \text{and } y(1) = 0 \text{ if } p = 2.$

So we see that in particular that $y$ is $\mathbb{F}$-linear. Let

$$\rho = \begin{pmatrix} \alpha & 0 \\ 0 & \alpha^{-1} \end{pmatrix},$$

where $\alpha$ is as in the proposition. The action of $H/Q$ on cocycles is given as follows: the cocycle $\alpha\xi$ maps $\sigma \mapsto \rho^{-1}\xi(\rho\sigma\rho^{-1})$. An easy computation now shows that a cocycle class $[\xi]$ represented by a cocycle $\xi$ is invariant under $H/Q$ when there is a coboundary $\eta$, given by $(m_1, m_2) \in V$, such that for all $\sigma_\lambda \in Q$,
$$\rho^{-1}\xi(\sigma_{\alpha^2\lambda}) = \xi(\sigma_\lambda) + \eta(\sigma_\lambda).$$
Let now $\tilde{\eta}$ be the coboundary given by $(0, \frac{m_2}{\alpha-1}) \in V$. An easy computation shows that if we replace $\xi$ by $\xi - \tilde{\eta}$, then for this $\xi$ the following equation holds:
$$(2) \quad \rho^{-1}\xi(\sigma_{\alpha^2\lambda}) = \xi(\sigma_\lambda).$$
This $\xi$ represents the class $[\xi]$ uniquely and the relations read in coordinates:
$$(2a) \quad \alpha^{-1}x(\alpha^2\lambda) = x(\lambda)$$
$$(2b) \quad \alpha y(\alpha^2\lambda) = y(\lambda).$$
So if we let $W$ be the $\mathbb{F}$-vectorspace consisting of tuples $(x, y)$, where $x : \mathbb{F} \longrightarrow \mathbb{F}$ is subject to the relations $(1a)$ and $(2a)$ and $y : \mathbb{F} \longrightarrow \mathbb{F}$ is subject to the relations $(1b)$ and $(2b)$, then $\dim_\mathbb{F} H^1(H, V) = \dim_\mathbb{F} W$. We will compute the latter dimension.

$(1b)$ and $(2b)$ imply that $\alpha^3 y(1) = y(1)$. So either $y = 0$ or $\alpha$ has order 3 and then $y$ is determined by $y(1)$.

If we let $\lambda_1, \ldots, \lambda_l$ be generators of $\mathbb{F}$ over $\mathbb{F}_p[\alpha]$, where $l = \dim_{\mathbb{F}_p[\alpha]} \mathbb{F}$, then by $(1a)$ and $(2a)$, $x$ is determined by $x(\lambda_i)$, with $i = 1, \ldots, l$. Hence $\dim_\mathbb{F} W \leq l + \delta$.

Suppose that $y(1) = 0$, then $x$ is $\mathbb{F}_p$-linear. Let $h$ be the minimal polynomial of $\alpha^2$, then for each $\lambda \in \mathbb{F}$, we have $0 = x(h(\alpha^2)\lambda) = h(\alpha)x(\lambda)$. Hence if $h(\alpha) \neq 0$, i.e. if $h$ is not the minimal polynomial of $\alpha$, then $x(\lambda) = 0$. Note that $h$ is the minimal polynomial of $\alpha$ iff $\alpha^{p^j} = \alpha^2$ for some $j \in \mathbb{N}$. Moreover if $\delta = 1$, then the order of $\alpha$ is 3, i.e. $\alpha^2 + \alpha + 1 = 0$ and in particular this means that $\alpha^p = \alpha^2$. So we see that if $\alpha^{p^j} \neq \alpha$, then $x(\lambda) = y(\lambda) = 0$, hence $\dim_\mathbb{F} H^1(H, V) = 0$.

Suppose now that $\alpha^{p^j} = \alpha^2$ and let $y(1) = 0$, then $x$ is not only $\mathbb{F}_p$-linear, but even $\mathbb{F}_p[\alpha]$-semi linear. This means that the $\mathbb{F}$-subspace of $W$ consisting of the tuples $(x, y)$, with $y = 0$, has dimension $\dim_{\mathbb{F}_p[\alpha]} \mathbb{F}' = l$.

If $\delta = 0$, then this subspace equals $W$ and we see $\dim_\mathbb{F} H^1(H, V) = l + \delta$. If $\delta = 1$, then the dimension of $W$ is either $l$ or $l + 1$. So let $\delta = 1$, then $p > 2$ and $\text{ord}(\alpha) = 3$. Suppose that $y(1) \neq 0$ and let $x : \mathbb{F} \longrightarrow \mathbb{F}$ be given by $x(\lambda) = c\lambda^2$, where $c = \frac{1}{2}y(1)$. Then one checks easily that $x$ has property $(1a)$. And because $\text{ord}(\alpha) = 3$, it has property $(2a)$ as well. This shows that $W$ contains an element $(x, y)$ with $y \neq 0$, hence $\dim_\mathbb{F} H^1(H, V) = l + \delta$. □

In the following lemma, we show that $\ker(\text{Res}_{\langle \sigma \rangle})$ does not depend on the choice of $1 \neq \sigma \in Q$.

**Lemma 4.5.** *For all $\sigma, \tau \in Q$, with both $\sigma$ and $\tau$ not the identity, we have $\ker(\text{Res}_{\langle \sigma \rangle}) = \ker(\text{Res}_{\langle \tau \rangle})$.*

*Proof.* For $p = 2$, by proposition (4.3) $H^1(\langle\sigma\rangle, V) = 0$, hence $\ker(\mathrm{Res}_{\langle\sigma\rangle}) = H^1(H, V)$ for all $\sigma \in Q$.

Let now $p > 2$. Note that
$$\ker(\mathrm{Res}_{\langle\sigma\rangle}) = H^1(H, V) \cap \ker(H^1(Q, V) \longrightarrow H^1(\langle\sigma\rangle, V)),$$
because $H^1(H, V) \hookrightarrow H^1(Q, V)$, so we may assume that $H = Q$.

Clearly, if $\sigma$ and $\tau$ are linearly dependent over $\mathbb{F}_p$, then $\langle\sigma\rangle = \langle\tau\rangle$, so $\ker(\mathrm{Res}_{\langle\sigma\rangle}) = \ker(\mathrm{Res}_{\langle\tau\rangle})$. If $\sigma$ and $\tau$ are independent over $\mathbb{F}_p$, then we may extend them to a basis $\langle\sigma_1, \ldots, \sigma_k\rangle$, with $\sigma = \sigma_1, \tau = \sigma_2$ and $k = \dim_{\mathbb{F}_p} Q$. We write $V = \mathbb{F}e_1 + \mathbb{F}e_2$ such that the $\sigma_i$'s are upper triangular on the basis $\{e_1, e_2\}$. Note that the kernel of $\mathrm{Res}_{\langle\sigma_i\rangle}$ is the image of
$$H^1(Q/\langle\sigma_i\rangle, V^{\langle\sigma_i\rangle}) \simeq \mathrm{Hom}(Q/\langle\sigma_i\rangle, \mathbb{F} \cdot e_1)$$
under the injective inflation map. The inflation map
$$\mathrm{Inf} : \mathrm{Hom}(Q/\langle\sigma_i\rangle, \mathbb{F} \cdot e_1) \longrightarrow H^1(Q, V)$$
is given by $\hat{\xi} \mapsto [\xi]$, such that if $\hat{\xi}([\sigma_j]) = a_j \in \mathbb{F}, j \neq i$, then $\xi$ is the cocycle given by
$$\xi(\sigma_j) = (a_j, 0) \in V, \quad j \neq i \quad \text{and} \quad \xi(\sigma_i) = (0, 0).$$
Now we will show that $\ker(\mathrm{Res}_{\sigma_k}) \subset \ker(\mathrm{Res}_{\sigma_l})$. If $[\xi] \in \ker(\mathrm{Res}_{\sigma_k})$, it comes from a $\hat{\xi}$ as mentioned above. Now let $\eta$ be a coboundary given by $(m_1, m_2) \in V$, hence
$$\eta : \sigma_i \mapsto (\lambda_i m_2, 0) \quad \text{for all} \quad i.$$
We choose $m_2$ such that $\lambda_l m_2 + a_l = 0$, then by construction there is a $\tilde{\xi} \in \mathrm{Hom}(Q/\langle\sigma_l\rangle, \mathbb{F} \cdot e_1)$ such that $\mathrm{Res}_{\sigma_l}(\tilde{\xi}) = [\xi + \eta]$. This shows that $[\xi] = [\xi + \eta] \in \ker(\mathrm{Res}_{\sigma_l})$. □

Recall that for any place $\omega$ of $K_0$ the map $\mathrm{Res}_\omega$ is the restriction map
$$\mathrm{Res}_\omega : H^1(G_{K_0}, V) \longrightarrow H^1(D_\omega, V).$$

**Proposition 4.6.** *Let $\varphi : A \longrightarrow K\{\tau\}$ be a Drinfeld module of rank 2 and let $H = G_{K_0} \cap \mathrm{Sl}_2(\mathbb{F})$ be of type (6), with $p$-Sylow group $Q$. Let $1 \neq \sigma \in Q$, then*
$$S(a, K) = \ker(H^1(G_{K_0}, V) \longrightarrow H^1(\langle\sigma\rangle, V)) \bigcap_{\omega : p^2 | \#D_\omega} \ker(\mathrm{Res}_\omega),$$
*where the intersection is taken over places $\omega$ of $K_0$. This intersection is finite.*

*Proof.* Suppose that $\omega$ is any place of $K_0$. If $p \nmid \#D_\omega$, then $\ker(\mathrm{Res}_\omega) = H^1(G_{K_0}, V)$, because $H^1(D_\omega, V) = 0$.

If $p \mid \#D_\omega$, but $p^2 \nmid \#D_\omega$, then
$$\ker(\mathrm{Res}_\omega) \subset \ker(H^1(G_{K_0}, V) \longrightarrow H^1(\langle\sigma\rangle, V)),$$

because $H^1(D_\omega, H) \hookrightarrow H^1(\langle \sigma \rangle, V)$. By Chebotarev's density theorem, it follows that there exists a place $\omega$ of $K_0$ with $D_\omega \cong \langle \sigma \rangle$, from which the theorem follows.

For the finiteness of the intersection, note if $p^2 \mid \#D_\omega$, then $\omega$ is ramified and there are only finitely many ramified places. $\square$

*Remark* 4.7. Clearly, $\ker(H^1(G_{K_0}, V) \longrightarrow H^1(\langle \sigma \rangle, V)) \subset \ker(\text{Res}_{\langle \sigma \rangle})$, with $\text{Res}_{\langle \sigma \rangle} : H^1(H, V) \longrightarrow H^1(\langle \sigma \rangle, V)$ as before. Hence (4.6) combined with (4.4) and (4.3) gives a bound on $\dim_{\mathbb{F}} S(a, K)$.

**Corollary 4.8.** *If $\varphi$ is a Drinfeld module of rank 2 over $\mathbb{F}_p$, $p > 2$ prime and $(a) \subset A$ is a prime ideal of degree 1, then $S(a, K) = 0$.*

*Proof.* In this case we have by theorem (4.6) that $S(a, K) \cong \text{Res}_{\langle \sigma \rangle}$ and because $H^1(G_{K_0}, V)$ embeds into $H^1(Q, V)$, we have by lemma (4.5), that $\dim_{\mathbb{F}} S(a, K) \leq -1 + \dim_{\mathbb{F}_p} Q$. Note that now $G_{K_0} \subset \text{Gl}_2(\mathbb{F}_p)$, hence $\dim_{\mathbb{F}_p} Q = 1$. $\square$

The proof of theorem (6.1) is similar to the proof of this corollary.

## 5. Examples of non-trivial $S(a, K)$

**Example 5.1.** Let $A = \mathbb{F}_q[t]$, where $\text{char}(\mathbb{F}_q) > 2$ and let $K = \mathbb{F}_q(t)$. Let $\varphi : A \longrightarrow K\{\tau\}$ be a Drinfeld module of rank 2 given by $\varphi_t = t + t\tau + t^2\tau^2$. We let $K_0 = K(\ker(\varphi_t))$, then $\text{Gal}(K_0/K) \subset \text{Gl}_2(\mathbb{F}_q)$.

We consider the decomposition group $D_t$. Clearly, the Newton polygon of $\varphi_t(Z)$ has two slopes, namely 0 and $\frac{1}{q(q-1)}$. To factor $\varphi_t(Z)$ in $\mathbb{F}_q((t))[Z]$, we need at least a completely ramified extension of degree $q(q-1)$ of $\mathbb{F}_q((t))$. Hence $D_t \cap \text{Sl}_2(\mathbb{F}_q)$ contains a subgroup of $q$ elements. If we compare this with the classification of subgroups of $\text{Sl}_2(\mathbb{F}_q)$ in section (4), we see that $D_t$ contains $\mathbb{F}_q$ as a subgroup. We conclude that $\text{Gal}(K_0/K)$ contains a subgroup isomorphic to $\begin{pmatrix} 1 & \mathbb{F}_q \\ 0 & 1 \end{pmatrix}$.

Recall that a system of Artin-Schreier equations over some $\mathbb{F}_p$-field $M$

$$\begin{cases} z_1^p - z_1 = f_1 & f_1 \in M \\ \quad \vdots \\ z_n^p - z_n = f_n & f_n \in M, \end{cases}$$

is called *independent* over $M$, if for all $\lambda_i \in \mathbb{F}_p$, such that not all $\lambda_i$ are 0, the equation $z^p - z = \sum_{i=1}^n \lambda_i f_i$ has no solutions in $M$. Such a system gives rise to a tower of field extensions $M = M_0 \subset M_1 \subset \ldots \subset M_n$ where the extension $M_i/M_{i-1}$ is given by $z_i^p - z_i = f_i$ and is of degree $p$.

**Example 5.2.** Let $q = p^k$, $p > 2$ and let

$$\varphi : \mathbb{F}_q[t] \longrightarrow K\{\tau\},$$

be a Drinfeld module of rank 2, such that $\mathrm{Gal}(K(\ker(\varphi_t)/K)$ contains a subgroup isomorphic to

$$H = \begin{pmatrix} 1 & \mathbb{F}_q \\ 0 & 1 \end{pmatrix},$$

e.g. the examples of (5.1). We write

$$\varphi_t = t + tc_1\tau + tc_2\tau^2, \quad c_1, c_2 \in K,$$

hence $K_0$ is the splitting field of the equation

$$(1) \quad 1 + c_1 Z^{q-1} + c_2 Z^{q^2-1} = 0.$$

Furthermore there exist elements $P, Q \in K_0$, such that $\ker(\varphi_t) = \mathbb{F}_q \cdot P + \mathbb{F}_q \cdot Q$. We let $L = K_0^H$. Because $\mathrm{Gal}(K_0/L) = H$, we may assume that $P \in L$.

If we substitute $U = Z^{q-1}$ in (1), we get

$$(2) \quad 1 + c_1 U + c_2 U^{q+1} = 0.$$

Let $L_1$ be the splitting field of (2), then we have the field inclusion $L \subset L_1 \subset K_0$, where the latter field extension is given by the equation $U = Z^{q-1}$. This implies that $[K_0 : L_1] \mid q-1$, but $[K_0 : L] = \#H = q$, hence $L_1 = K_0$. This shows that $K_0$ is the splitting field of (2) over $L$.

Because $P \in L$, we already know a solution of (2), namely $u = P^{q-1}$. Substituting $V = U - u$ in (1) gives

$$1 + c_1(V + u) + c_2(V + u)(V^q + u^q) = c_1 V + c_2 V u^q + c_2 u V^q + c_2 V^{q+1}.$$

Subsequently we divide out $V$ and substitute $W = V^{-1}$, which shows that $K_0/L$ is the splitting field of the equation

$$(3) \quad W^q + \frac{c_2 u}{c_1 + c_2 u^q} W + \frac{c_2}{c_1 + c_2 u^q} = 0.$$

To simplify this equation a little more, we consider it over the extension $L(b)$ of $L$, with $b^{q-1} = -\frac{c_2 u}{c_1 + c_2 u^q}$. Because $[L(b) : L] \mid q-1$, the degree of $M := L(b)$ over $L$ is relatively prime to $q$, hence the splitting field $M_0 := K_0(b)$ of (3) over $M$ also has Galois group $\mathrm{Gal}(M_0/M) \cong H$.

Substituting $bX = W$ in (3) gives

$$X^q - X = f, \quad \text{where} \quad f = \frac{1}{bu}.$$

The following theorem shows that $S(a, K)$ can be arbitrarily large.

**Theorem 5.3.** *For any $k \in \mathbb{N}_{>0}$, there exists a function field $K$, a Drinfeld module $\varphi : A \longrightarrow K\{\tau\}$ and a prime ideal $(a) \subset A$, such that $\dim_\mathbb{F} S(a, K) = k$.*

*Proof.* Let $q = p^k$ for some integer $k > 1$ and $p > 2$ a prime. The computations of example (5.1) and (5.2) show that there is a Drinfeld module $\varphi$ over some function field $M$, such that $M_0 = M(\ker(\varphi_t))$ is a Galois extension with Galois group $H = \begin{pmatrix} 1 & \mathbb{F}_q \\ 0 & 1 \end{pmatrix}$ and moreover this extension $M_0/M$ is an Artin-Schreier extension given by

$$(4) \qquad X^q - X = f, \quad f \in M.$$

This extension is also given by the system of Artin-Schreier equations

$$(5) \qquad \begin{cases} x_1^p - x_1 = \beta_1 f \\ \qquad \vdots \\ x_k^p - x_k = \beta_k f \end{cases},$$

where the $\beta_i \in \mathbb{F}_q$ are linearly independent over $\mathbb{F}_p$. To see this, write $z = \sum_{i=1}^k \alpha_i x_i$, with $\alpha_i \in \mathbb{F}_q$. An easy computation shows that $z$ is a solution of $X^q - X = f$ if and only if

$$\begin{pmatrix} \beta_1 & \cdots & \beta_k \\ \beta_1^p & \cdots & \beta_k^p \\ \vdots & & \vdots \\ \beta_1^{p^{k-1}} & \cdots & \beta_k^{p^{k-1}} \end{pmatrix} \begin{pmatrix} \alpha_1 \\ \alpha_2 \\ \vdots \\ \alpha_k \end{pmatrix} = \begin{pmatrix} 1 \\ 0 \\ \vdots \\ 0 \end{pmatrix}.$$

Because this matrix is invertible (cf. [Gos96]), it follows that (5) is indeed equivalent to $X^q - X = f$.

Now we consider the extension $M(z_1)/M$, given by $z_1^p - z_1 = \beta_1 f - g_1$, where we choose $g_1$ as follows: for the finitely many places $\nu$ of $M$, for which $v_\nu(f) < 0$, we let $v_\nu(g_1) > 0$ and for one place $\nu_0$ for which $v_{\nu_0}(f) > 0$, we let $v_{\nu_0}(g_1) = -1$. Such a $g_1$ exists, cf. [Bou89], corollary VI.2.1. The condition $v_{\nu_0}(g) = -1$, makes sure that system given by (5) and the equation for $z_1$ is independent: Namely, by Hensel's lemma all equations of (5) have solutions in $M_{\nu_0}$, whereas the equation for $z_1$ gives rise to a totally ramified extension of degree $p$.

We conclude that (5) is still independent over $M(z_1)$. Now we consider $M(z_1, z_2)/M(z_1)$ given by $z_2^p - z_2 = \beta_2 f + g_2$. We choose $g_2$ in the same way as we chose $g_1$, with $M$ replaced by $M(z_1)$. This implies that (5) is independent over $M(z_1, z_2)$. By repeating this we get that (5) is independent over the field $M(z_1, \ldots, z_{k-1})$.

Let $L/M(z_1, \ldots, z_{k-1})$ be the field extension given by (5), then its Galois group is $H$. Let $\nu$ be a place of $M(z_1, \ldots, z_{k-1})$ and let $v_\nu$ be its corresponding valuation. Let $\omega$ be a place of $L$ lying above $\nu$. We distinguish the following cases.

(a) $v_\nu(f) > 0$, in this case we see that the equations of (5) are over the residue field given by $x_i^p - x_i = 0$, hence they split completely over the residue field. Hensel's lemma implies that $D_\omega$ is trivial.

(b) $v_\nu(f) = 0$, then also $v_\nu(\beta_i f) = 0$, hence all equations of (5) are over the residue field given by $x_i^p - x_i = \alpha_i$, with $\alpha_i$ in the residue field. Hence all equations only give rise to a residue field extension. This shows that $\nu$ is in $L$ and thus $D_\omega$ is cyclic and can have at most $p$ elements, because the elements of $H$ have at most order $p$.

(c) $v_\nu(f) < 0$. Note that the equations $x_i^p - x_i = \beta_i f$ are equivalent to $y_i^p - y_i = g_i$ by substituting $y_i = z_i - x_i$, for $i = 1, \ldots, k-1$. Because by construction $v_\nu(g_i) > 0$, it follows that these equations give a trivial extension at $\nu$. So only the equation $x_k^p - x_k = \beta_k f$ can give rise to a non-trivial extension, but this extension has at most degree $p$, hence $D_\omega$ can at most have $p$ elements.

We see that at any place $\omega$, the decomposition group $D_\omega$ has at most $p$ elements. This means that for the non-trivial $D_\omega$, the kernel $\ker(\mathrm{Res}_\omega)$ has dimension

$$\dim_{\mathbb{F}_q} \ker(\mathrm{Res}_\omega) = -1 + \dim_{\mathbb{F}_q} H^1(H,V) = k-1,$$

by proposition (4.3). Hence it follows by proposition (4.5), that

$$\dim_{\mathbb{F}_q} S(t, M(z_1, \ldots, z_{k-1})) = k-1.$$

$\square$

## 6. THE ELLIPTIC CURVE CASE

In this section we will treat the analogous problem for elliptic curves. Let $E$ be an elliptic curve over some number field $K$ and let $p \in \mathbb{N}$ be a prime number. For any $P \in E(K)$ we denote $K_P = K(p^{-1}P)$. In this section we will prove the following theorem:

**Theorem 6.1.** *Let $E$ be an elliptic curve over a number field $K$, let $p$ be a prime number, then the kernel*

$$S(p, K) = \ker \left( E(K)/pE(K) \longrightarrow \prod_\nu E(K_\nu)/pE(K_\nu) \right),$$

*where $\nu$ runs through the places of $K$, is trivial.*

As before, we will write $G_{K_0} = \mathrm{Gal}(K_0/K)$. Because $E[p](\overline{K}) \cong \mathbb{F}_p P + \mathbb{F}_p Q$, with $P, Q \in K_0$, we see that

$$G_{K_0} \hookrightarrow \mathrm{Gl}_2(\mathbb{F}_p).$$

Clearly $K_0 \subset K_P$. We will denote $V = E[p](K_0) = E[p](K_{0,\omega})$.

**Proposition 6.2.** *For every $P \mod pE(K) \in S(p,K)$, we have $K_P = K_0$. In particular $S(p, K_0) = 0$.*

*Proof.* For every $Q \in pE(K)$, clearly $K_Q = K_0$, hence $K_P$ only depends on the class $[P] = P + pE(K)$. Furthermore, if $P \in S(p,K)$, then $P \in pE(K_\nu)$ for every place $\nu$ of $K$. If we let $\nu_0$ be a place of $K_0$ lying above $\nu$ and $\nu_P$ be a place of $K_P$ lying above $\nu_0$, then this implies that $(K_\nu)_P \subset K_{0,\nu_0}$. This gives rise to an embedding $K_0 \subset K_P \subset K_{0,\nu_0}$, hence $\deg(\nu_P/\nu_0) = 1$. By lemma (2.2) it follows that $K_P = K_0$.

Now $S(p, K_0) = 0$ as in (2.3). $\square$

**Proposition 6.3.** *We have that*
$$S(p, K) \subset \bigcap_\omega \ker(\mathrm{Res}_\omega),$$
*where the intersection is taken over all places $\omega$ of $K_0$ and the map $\mathrm{Res}_\omega$ is the restriction map*
$$\mathrm{Res}_\omega : H^1(G_{K_0}, V) \longrightarrow H^1(D_\omega, V),$$
*with $D_\omega$ the decomposition group at $\omega$.*

*Proof.* As in the proof of (2.4), we have the following diagram with exact rows and columns:

$$\begin{array}{ccccccc}
 & & 0 & & 0 & & \\
 & & \downarrow & & \downarrow & & \\
0 & \longrightarrow C \longrightarrow & S(p, K) & \longrightarrow & S(p, K_0) = 0 & & \\
 & & \downarrow & & \downarrow & & \\
0 & \longrightarrow \Phi \longrightarrow & E(K)/pE(K) & \longrightarrow & E(K_0)/pE(K_0) & & \\
 & & \downarrow & & \downarrow & & \\
0 & \longrightarrow B \longrightarrow & \prod_\nu E(K_\nu)/pE(K_\nu) & \longrightarrow & \prod_\omega E(K_{0,\omega})/pE(K_{0,\omega}). & &
\end{array}$$

Applying the same arguments of Galois cohomology as in the proof of proposition (2.4) we obtain injections $E(L)/pE(L) \hookrightarrow H^1(\mathrm{Gal}(\overline{K}/L), E[p](\overline{K})$ for $L = K$ and $L = K_0$. This gives rise to an embedding $\Phi \hookrightarrow H^1(G_{K_0}, V)$. Arguing in the same way we get an embedding
$$B \hookrightarrow \prod_\nu H^1(D_\omega, E[p](K_{0,\omega}),$$
where the product runs over all places $\nu$ of $K$. This implies that
$$C \hookrightarrow \bigcap_\omega \ker(\mathrm{Res}_\omega).$$
□

Let $H = G_{K_0} \cap \mathrm{Sl}_2(\mathbb{F}_p)$. As before we have
$$H^1(G_{K_0}, V) \hookrightarrow H^1(H, V).$$
Because $H \subset \mathrm{Sl}_2(\mathbb{F}_p)$ we have that $H$ is one of the following, cf. the classification in section (4):

(1) $p \nmid \#H$.
(2) $D_2$; in this case $p = 2$.
(3) $\mathrm{Sl}_2(\mathbb{F}_p)$.
(4) $H$ is a Borel group, i.e. $H$ has a cyclic normal subgroup $Q = \langle \sigma \rangle$ of order $p$ and $H/Q$ is cyclic of order dividing $p - 1$.

**Proposition 6.4.** *If $H$ is of type (1), (2) or (3), then $H^1(G_{K_0}, V) = 0$.*

*Proof.* Except for the case $H = \mathrm{Sl}_2(\mathbb{F}_2)$, this follows from (4.1). So let $H = \mathrm{Sl}_2(\mathbb{F}_2)$ and let $\sigma \in H$ be an element of order 2. The group $H^1(\langle\sigma\rangle, V) = 0$ - this is (4.3), with $p = 2$ and $Q = \langle\sigma\rangle$. Consider the restriction-corestriction sequence

$$H^1(H, V) \xrightarrow{\mathrm{Res}} H^1(\langle\sigma\rangle, V) \xrightarrow{\mathrm{Cor}} H^1(H, V).$$

Then $\mathrm{Cor} \circ \mathrm{Res} = [H : \langle\sigma\rangle]$, hence

$$[H : \langle\sigma\rangle] H^1(H, V) = 0.$$

Because $[H : \langle\sigma\rangle]$ is relatively prime to 2, it follows that $H^1(H, V) = 0$. □

*Proof. Of theorem (6.1).* Because $S(a, K) \hookrightarrow H^1(H, V)$, we know by (6.4) that if $H$ is not of type (4), then $S(a, K) = 0$. Suppose now that $H$ is of type (4) and let $\sigma \in H$ be an element of order $p$. Then by Chebotarev and (6.3), we have $S(a, K) \hookrightarrow \ker(\mathrm{Res}_{\langle\sigma\rangle})$. By (4.3), this kernel has dimension $-1 + \dim_{\mathbb{F}_p} Q = -1 + 1 = 0$. □

Vakgroep Wiskunde RuG, P.O. Box 800, 9700 AV Groningen, the Netherlands